\documentclass [14pt]{amsart}
\usepackage{amsmath}
\usepackage{amssymb}
\usepackage{mathrsfs} 
\usepackage{enumitem}
\usepackage{amscd}
\usepackage{epsfig}
\usepackage{amsxtra}
\usepackage{pifont}
\usepackage{mathabx}
\usepackage{MnSymbol}
\usepackage{accents}
\usepackage{scalerel,stackengine}
\usepackage{xcolor}
\usepackage{extarrows}
\stackMath

\usepackage[
backend=biber, 
hyperref=true, 
style=ext-alphabetic, 
citexref=true,
url=false,
doi=false,
isbn=false,
eprint=false,
articlein=false,
giveninits=true,
date=year,
maxnames=5,
minnames=5,
maxalphanames=5,
minalphanames=5
]{biblatex}  %% Main settings

\usepackage{hyperref}
\hypersetup{
	colorlinks=true,
	allcolors=blue
}

\DeclareFieldFormat[online]{date}{Last Modified: #1}

\DeclareFieldFormat[article]{pages}{#1}

\DeclareFieldFormat
[article,inbook,incollection,inproceedings,patent,thesis,unpublished]
{title}{#1}

\DeclareFieldFormat
[article,inbook,incollection,inproceedings,patent,thesis,unpublished]
{titlecase:title}{\MakeSentenceCase*{#1}}

\DeclareFieldFormat[article,periodical]{number}{\mkbibparens{#1}}

\DeclareFieldFormat{issuedate}{#1}

\newcommand\reallywidehat[1]{%
\savestack{\tmpbox}{\stretchto{%
  \scaleto{%
    \scalerel*[\widthof{\ensuremath{#1}}]{\kern-.6pt\bigwedge\kern-.6pt}%
    {\rule[-\textheight/2]{1ex}{\textheight}}%WIDTH-LIMITED BIG WEDGE
  }{\textheight}%
}{0.5ex}}%
\stackon[1pt]{#1}{\tmpbox}%
}

\newcommand\reallywidecheck[1]{%
\savestack{\tmpbox}{\stretchto{%
  \scaleto{%
    \scalerel*[\widthof{\ensuremath{#1}}]{\kern-.6pt\bigwedge\kern-.6pt}%
    {\rule[-\textheight/2]{1ex}{\textheight}}%WIDTH-LIMITED BIG WEDGE
  }{\textheight}%
}{0.5ex}}%
\stackon[1pt]{#1}{\scalebox{-1}{\tmpbox}}%
}

\allowdisplaybreaks
\newcommand{\supp}{\mbox{\rm supp}}

\newcommand{\RR}{{\mathbb R}}

\newcommand{\ZZ}{{\mathbb Z}}
\newcommand{\CC}{{\mathbb C}}

\newcommand{\NN}{\mathbb N}

\newcommand{\cM}{{\mathcal M}}

\newcommand{\cS}{{\mathcal S}}

\newcommand{\dd}{\mbox{d}}

\newcommand{\Cc}{C_{\mathsf{c}}}
\newcommand{\Ccinf}{C_{\mathsf{c}}^\infty}

\newcommand{\STM}{\mathcal S \mathcal M}

\def\hatf{\reallywidehat{f}}

\def\hatmu{\reallywidehat{\mu}}

\def\hatgam{\reallywidehat{\gamma}}

\newcommand{\seq}[1]{\left\{#1\right\}_{n\in\NN}}

\newcommand{\card}{{\rm card}}
\def\lt{\left}
\def\rt{\right}

\def\intd{\int_{\RR^d}}

\def\deq{:=}

\numberwithin{equation}{section}

\newtheorem{theorem}{Theorem}[section]
\newtheorem{lemma}[theorem]{Lemma}
\newtheorem{proposition}[theorem]{Proposition}
\newtheorem{corollary}[theorem]{Corollary}
\newtheorem{fact}[theorem]{Fact}

\theoremstyle{definition}

\newtheorem{definition}[theorem]{Definition}
\newtheorem{example}[theorem]{Example}
\newtheorem{remark}[theorem]{Remark}

\begin{document}
\title{Circles in diffraction}

\author{Emily R. Korfanty}
\address{Department of Mathematical and Statistical Sciences, University of Alberta, Edmonton, AB,  T6G 2G1, Canada \\
}
\email{ekorfant@ualberta.ca}
\urladdr{}

\author{Nicolae Strungaru}
\address{Department of Mathematical Sciences, MacEwan University \\
10700 -- 104 Avenue, Edmonton, AB, T5J 4S2, Canada\\
and \\
Institute of Mathematics ``Simon Stoilow''\\
Bucharest, Romania}
\email{strungarun@macewan.ca}
\urladdr{https://sites.google.com/macewan.ca/nicolae-strungaru/home}
\maketitle
\begin{abstract}
Given a Fourier transformable measure in two dimensions, we find a formula for the intensity of its Fourier transform along circles. In particular, we obtain a formula for the diffraction measure along a circle in terms of the autocorrelation measure. We look at some applications of this formula.
\end{abstract}

\section{Introduction}

Since the early 1900s, crystallographers have used diffraction experiments to determine the structure of crystals. Until the discovery of quasicrystals \cite{SBGC1984} by Dan Shechtman in 1984, pure point diffraction was taken as evidence of an underlying periodic structure in the material.  Quasicrystals, however, exhibit sharp diffraction patterns, and typically show symmetries that are impossible in materials with a periodic structure.  Shechtman was later awarded the 2011 Nobel Prize in Chemistry for the discovery, which also inspired a new interest in aperiodic structures and mathematical diffraction theory. For an overview of this area, we recommend the monographs \cite{BG2013,BG2017}.  

The rigorous mathematical treatment of diffraction was introduced by Hof in \cite{Hof1995}. He defined the diffraction measure as the Fourier transform of the so-called \emph{autocorrelation measure} of the underlying structure. 
%Over the last few years diffraction has been well understood for systems supp
%Since then, diffraction has been well understood when the underlying 
%point set $\Lambda$ has finite local complexity \cite{BG2013} (i.e. when $\Lambda - \Lambda$ is closed and discrete)
%tiling has finite local complexity,  especially in the case of Meyer set support \cite{Str2005,Str2013,Str2014,Str2017, Str2020}. However, the techniques used in these papers cannot be applied for tilings with infinite local complexity.
Diffraction is now well understood for systems with pure point spectrum \cite{LSS2020,LSS2024}, and systems coming from cut-and-project schemes \cite{Sch2000,BM2004,RS2017,Str2005,Str2017, Str2020, BHS2017, KR2016}. 

On the other hand, the diffraction for systems without finite local complexity (FLC) is, in general, much less understood. One such example, the pinwheel tiling \cite{Rad1994}, is produced by a simple 2-dimensional substitution. The tiles in the pinwheel tiling appear in infinitely many orientations, and therefore the pinwheel does not have FLC.  Its diffraction is rotationally invariant  \cite{MPS2006}, and has a Bragg peak of unit intensity at the origin, but nothing else is known about its diffraction. In particular, it is not known if the diffraction is absolutely continuous, singular continuous or mixed. 

There are many examples of rotational invariant diffraction patterns. For example, any substitution with \textit{statistical circular symmetry}, meaning that the orientations of the tiles distribute uniformly around the unit circle (\cite{Fre2008} \cite{Sad1998}, \cite{FGH}), leading to a rotational invariant diffraction pattern. More generally, the diffraction of any rotational invariant measure is also rotational invariant. 

Rotational invariance also appears in powder diffraction, which is a slightly different diffraction technique, where the sample is ground into many smaller samples prior to the usual diffraction process.  During the grinding process, the samples tend to rotate in random directions. This random distribution of rotations on the unit circle is typically statistically circular, and hence the powder diffraction is rotational invariant.

%. This is because the pinwheel tiling exhibits \textit{statistical circular symmetry}, meaning that the orientations of the tiles distribute uniformly around the unit circle. 
%Aside from this, there is little else known about the pinwheel diffraction, except that there is a Bragg peak of unit intensity at the origin.  

%More is known about the autocorrelation of $\Lambda_{\pin}$; indeed, we know that the autocorrelation consists only of circles \cite{MPS2006}, and that the radii of these circles must be contained in the following countable set \cite{BFG2007a,BFG2007}:
%\[
%\{\text{radii of circles in }\gamma_\pin\}= \lt\{\|x-y\|\,:\, x,y\in\Lambda \rt\} \subseteq \lt\{\sqrt{\frac{p^2+q^2}{5^\ell}}\,:\, p,q,\ell \in \NN_0\rt\}\,.
%\]
%In theory, one can even determine the radii precisely up to any fixed $R>0$ by using a collared version of the pinwheel tiling \cite{FWW2014} and the renormalization procedure introduced in \red{\cite{?}}.

In \cite{BFG2007,BFG2007a,GD2011}, numerical approximations of the pinwheel diffraction were compared to a powder diffraction from $\ZZ^2$.  These numerical approximations suggest that the pinwheel diffraction may include bright circles with the same radii as the circles in the powder diffraction from $\ZZ^2$, which would imply the existence of a singular continuous component of the spectrum. Despite this observed similarity to the powder diffraction, it is still unknown whether or not there exists any circle in the pinwheel diffraction.  Besides numerical simulations, which are not too precise for this particular system,  we are not aware of any progress made toward understanding the diffraction of the pinwheel tiling in the last 20 years.

 %. Indeed, no such proof has been found, and the error in these approximations is potentially too large to form a confident conjecture.  To the authors' knowledge, there has been no significant progress in this direction in the last decade. 

Motivated by the above, in this paper we develop a method for computing the intensity of a circle in the diffraction. This may provide a new tool to study the diffraction of substitution tilings with statistical circular symmetry, such as the pinwheel and its related families \cite{Sad1998,Fre2008, FGH}.\textbf{}  In Section~\ref{sec:diffraction}, we derive the following explicit limit formula for the diffraction intensity $\hatgam(C_r)$ along the circle $C_r=\{ x \in \RR^2 :\|x \|=r\}$  in terms of the autocorrelation measure $\gamma$:
\begingroup
\def\thetheorem{\ref{thm:diffraction-formula}}
\begin{theorem}
    Let $\gamma$ be an autocorrelation measure on $\RR^2$ and let $r>0$. Then
    \begin{equation*}
        \hatgam(C_r) = \lim_{n\rightarrow\infty} \frac{2r\pi\sqrt{\pi}}{n}\int_{\RR^2} e^{-\pi^2\|y\|^2/n^2}J_0(2\pi r \|y\|) \,\dd\gamma(y)\,.
    \end{equation*} 
\end{theorem}
\addtocounter{theorem}{-1}
\endgroup
\noindent Here $J_0$ represents the Bessel function of the first kind (see Definition~\ref{def:J0}).

To obtain this result, our approach is as follows. Heuristically, we should have 
\[
\reallywidehat{\gamma}(C_r)= \reallywidehat{\gamma}(1_{C_r})\text{``=''} \gamma(\reallywidehat{1_{C_r}}) \,.
\]
Of course, the last equality above cannot hold, since $1_{C_r} =0$ almost everywhere with respect to the Lebesgue measure on $\RR^2$, and hence $\reallywidehat{1_{C_r}}=0$. To get around this issue, we instead approximate $1_{C_r}$ by a sequence of Schwartz functions of the form
\[
f_n:= g_n * \theta_r\,,
\]
where $\theta_r$ is the uniform probability measure on $C_r$ (see Definition~\ref{def:thetar}), and $g_n$ is an approximate identity sequence of Gaussian functions, renormalized such that $f_n$ converges to $1$ on $C_r$. Using the formula $\reallywidehat{\gamma}(f_n)=\gamma(\reallywidehat{f_n})$ and letting $n \to \infty$ we obtain Theorem~\ref{thm:diffraction-formula}. The Bessel function $J_0$ appears naturally in this approach through the Fourier transform of $\theta_r$; see Proposition~\ref{prop:circle-FT}.

Our paper is organized as follows. We start by reviewing the basic definitions and framework in  Section~\ref{sec:framework}. We use the Fourier transform of strongly tempered measures (see Definition~\ref{def:STM}), which extends the classical theory of Fourier transformable measures in $\RR^d$. In Section~\ref{sec:Bessel} we briefly review the Bessel functions $J_0$ and $I_0$ and some of their properties. As emphasized above, $J_0$ appears as the Fourier transform of $\theta_r$, while $I_0$ appears naturally in the convolution of $\theta_r$ with the Gaussian function(s).

%There are many tilings other than the pinwheel with rotationally invariant diffraction . The substitution rules that generate these tilings are more complicated than for the pinwheel tiling, and their diffractions are much less studied. In particular, it is also unknown for these tilings whether or not there exists a circle in the diffraction.  Thus, any progress towards better understanding the measure of a circle in diffraction warrants interest.

%Our results also provide formulae involving the Bessel function $J_0$ when applied for measures where the Fourier transform is known.  One particularly noteworthy example is when $\gamma = \delta_{\ZZ^2}$. In Section~\ref{sec:applications}, we prove the following about the sum of two squares function $r_2$, which has been a topic of great interest in number theory \red{ (should we provide some references here?)}:

%In Section~\ref{sec:basics}, we give some basic definitions and notations, and establish the general framework in which our main result Theorem~\ref{thm:diffraction-formula} holds.  More specifically, we introduce the space $\STM(\RR^d)$ of strongly tempered measures that have a distributional Fourier transform that is also a strongly tempered measure. Then, in Section~\ref{sec:Bessel}, we prove the specific facts about Bessel functions that we will need to prove Theorem~\ref{thm:diffraction-formula}.  

Section~\ref{sec:circles} provides the fundamental tool for proving Theorem~\ref{thm:diffraction-formula}.  We show, via a straightforward application of the dominated convergence theorem, that the intensity $\hatgam(C_r)$ can be calculated as the limit of integrals $\int_{\RR^2} f_n(x) \,\dd\hatgam(x)$, whenever when the sequence $f_n$ of Schwartz functions satisfies certain conditions (which we call approximate $r$-circles). Since $\int_{\RR^2} f_n(x) \,\dd\hatgam(x)=\int_{\RR^2} \reallywidehat{f_n}(y) \,\dd\gamma(y)$, we obtain a more general version of Proposition~\ref{prop:shelling} (see Corollary~\ref{cor:STM}).
%For this, we take $f_n$ to be a Schwartz function and use the fact that both $\gamma$ and $\hatgam$ are strongly tempered to ensure that $f_n$ and the Fourier transforms $\reallywidehat{f_n}$ are integrable with respect to $\hatgam$ and $\gamma$, respectively.  \red{Should we mention Feichtinger algebra?}  Then, we use the distributional Fourier transforms to conclude that
%$$
%\mu(C_r) = \lim_{n\rightarrow\infty} \int_{\RR^2} \reallywidehat{f_n}(y)\,\dd \mu(y)\,,
%$$
%for any $\mu\in\STM(\RR^2)$.
Theorem~\ref{thm:diffraction-formula} is then obtained by constructing in Section~\ref{sec:diffraction} an explicit approximate $r$-circle with easily calculable Fourier transform.

In Section~\ref{sec:applications}, we complete the paper by looking at several examples and applications of Theorem~\ref{thm:diffraction-formula}. In particular, the Poisson summation formula in $\RR^2$ yields the following result for the sum of two squares function:
\begingroup
\def\thetheorem{\ref{prop:shelling}}
\begin{proposition}
The sum of two squares function $r_2(n)$ satisfies the following:
   \begin{itemize}
       \item[(a)] For all $k \in \NN$ we have 
       \[
r_2(k)= \lim_{n\rightarrow\infty} \frac{2\pi \sqrt{\pi k}}{n} \sum_{m =0}^{\infty}r_2(m) e^{-\frac{\pi^2m}{n^2}}J_0(2\pi \sqrt{km})  \,.
       \]
       \item[(b)] For all $r>0$ such that $r^2 \notin \NN$ we have 
       \[
 \lim_{n\rightarrow\infty} \frac{2r \pi\sqrt{\pi}}{n} \sum_{m =0}^{\infty}r_2(m) e^{-\frac{\pi^2m\textbf{}}{n^2}}J_0(2\pi  r \sqrt{m})  =0 \,.
       \]
   \end{itemize}
\end{proposition}
\addtocounter{theorem}{-1}
\endgroup
In Appendix~\ref{sect:anotherapp}, we construct a different family of approximate $r$-circles via convolutions with annuli of small width.

\section{Preliminaries}\label{sec:basics}

\subsection{Definitions and notations}

Throughout this paper, $\RR^d$ denotes the Euclidean space of dimension $d\geq 1$. We will denote the standard Euclidean norm on $\RR^d$ by $\|\cdot\|$. We will primarily be interested in the case $d=2$. 

Given $R> 0$, we denote by $B_R$  the open ball 
\[
B_R \deq \{x\in\RR^d\,:\,\|x\| < R\}\,,
\]
of radius $R$, and by $\overline{B_R}$  the closed ball
\[
\overline{B_R}\deq \{x\in\RR^d\,:\,\|x\| \leq R\} \,.
\] 
We will use $1_A$ to denote the characteristic function of a set $A\subseteq \RR^d$, i.e.
\[
1_A(x) \deq \begin{cases}
1 & x \in A \\
0 & x \notin A\,.
\end{cases}
\]
\smallskip

As usual, $\Cc(\RR^d)$ denotes the space of continuous functions with compact support, and $\Ccinf(\RR^d)$ is the subspace of $\Cc(\RR^d)$ consisting of smooth (infinitely many times differentiable) functions.  The space of Schwartz functions is denoted by $\cS(\RR^d)$ and the space of tempered distributions is denoted by $\cS'(\RR^d)$.  

If $f$ is any bounded function, we will use the notation 
\[
\|f\|_\infty \deq \sup_{x\in\RR^d} |f(x)|\,,
\]
for the supremum norm of $f$.  

Let us recall that the \textit{Fourier transform} of a function $f\in L^1(\RR^d)$ is the function $\reallywidehat{f} : \RR^d \to \CC$ defined by
\[
\hatf(y) \deq \intd f(x) e^{-2\pi i x \cdot y} \,\dd x \qquad \forall y\in\RR^d\,.
\]
The map $f\mapsto\hatf$ is a bijection from $\cS(\RR^d)$ onto itself, and thus induces a mapping $\reallywidehat{\cdot}: \cS'(\RR^d) \to \cS'(\RR^d)$, via
 \begin{equation}\label{eq:TD-FT}
\reallywidehat{T}(f) \deq T(\reallywidehat{f}\,) \qquad \forall f\in\cS(\RR^d)\,.  
 \end{equation}
For a tempered distribution $T \in \cS'(\RR^d)$ we will refer to the tempered distribution $\reallywidehat{T}$ as the \emph{distributional Fourier transform of $T$}.

\subsection{Strongly tempered measures}\label{sec:framework}
In this paper, we prove our results for a space of measures with a certain property (see Definition~\ref{def:STM} below) that is compatible with diffraction theory.  Let us first recall the definition of a measure.

\begin{definition}\label{def:measure}
A \textit{Radon measure} on $\RR^d$ is a linear functional $\mu: \Cc(\RR^d)\rightarrow \CC$ with the property that for every compact set $K\subseteq \RR^d$, there exists a constant $a_K$ such that
\begin{equation}\label{eq:regular}
|\mu(f)|\leq a_K \|f\|_\infty\,,
\end{equation}
holds for every $f\in\Cc(\RR^d)$ with $\supp(f)\subseteq K$.  We write $\cM(\RR^d)$ for the space of all Radon measures on $\RR^d$. We will simply refer to Radon measures as measures.
\end{definition}

We will often use the so-called \textit{variation measure} $|\mu|$, where $\mu$ is some measure.  The definition of $|\mu|$ is provided in the following theorem:

\begin{theorem}[{{\cite[Lemma~C.2]{RS2017}}, \cite[\nopp 6.5.8]{Ped1995}}]\label{thm:variation}
For every measure $\mu$, there exists a smallest positive measure $|\mu|$ satisfying
\[
|\mu(f)|\leq |\mu|(|f|) \qquad \forall f\in\Cc(\RR^d)\,.
\]
Moreover, $|\mu|$ is the unique measure with the property that 
\[
|\mu|(f)= \sup \{|\mu(g)| : g \in \Cc(\RR^d), |g| \leq f \}\,,
\]
holds for all non-negative functions $f \in \Cc(\RR^d)$.
\qed
\end{theorem}

\medskip

Below, we will be using the Fourier transform in the distributional sense for measures. Let us first note  that, in general, $\int_{\RR^d} f(x) \,\dd \mu(x)$ may not make sense for a measure $\mu$ and Schwartz function $f$.  Since we will work with integrals of this type, we will need to restrict the class of measures we are studying. We start by recalling the following definition from \cite{BS2022}:

\begin{definition} 
%A measure $\mu\in\cM(\RR^d)$ is called \emph{tempered} if there exists a tempered distribution $T\in\cS'(\RR^d)$ such that 
%\[
%\int_{\RR^d} \varphi(x) \, \dd \mu(x) = T(\varphi) \qquad \forall \varphi \in \Cc^\infty(\RR^d) \,.
%\]
A measure $\mu$ is called\footnote{In \cite{AG1974} the authors call such a measure slowly increasing. The term strongly tempered is used in \cite{BS2022}.} \emph{strongly tempered} if  there exists a polynomial $P \in \CC[x_1, x_2, \ldots, x_d]$ such that
\[
\int_{\RR^d} \frac{1}{1+|P(x)|} \, \dd |\mu|(x) < \infty \,.
\]
\end{definition}
%\begin{remark}\cite[Theorem 2.7]{BS2022} Let $\mu$ be a positive measure.  Then $\mu$ is tempered if and only if $\mu$ is strongly tempered.
%\end{remark} 
The following fact about strongly tempered measures explains why this definition is important for us.
\begin{fact}[{\cite[Theorem 2.7]{BS2022}}]\label{fact:strongly-tempered}
A measure $\mu$ is strongly tempered if and only if we have $\int_{\RR^d} |f(t)| \dd |\mu|(t) < \infty$ for all $f\in\cS(\RR^d)$. Moreover, in this case, the map $T_\mu: \cS(\RR^d)\rightarrow\CC$ given by
\[
  T_\mu(f):=\int_{\RR^d} f(t) \dd \mu(t)\,,
\]
defines a tempered distribution. Furthermore, in this case, the measure $\mu$ and the tempered distribution $T_\mu$ uniquely determine each other.  
\end{fact}
Let us now discuss the Fourier transformability of strongly tempered measures in the distributional sense. Every strongly tempered measure has a distributional Fourier transform; however, this Fourier transform is a tempered distribution that is not necessarily a measure. Even if it is a measure, it may not be strongly tempered (see for example \cite[Prop.~7.1]{AG1974} or \cite[Sect.~3]{BS2022}).  Since our results rely on the distributional Fourier transform being a strongly tempered measure, we introduce the following definition. 

\begin{definition}\label{def:STM}
A measure $\mu$ is called \textit{Fourier transformable in the strongly tempered sense} if $\mu$ is a strongly tempered measure and there exists a strongly tempered measure $\nu$ such that, as tempered distributions, we have 
\[
\reallywidehat{T_{\mu}}= T_\nu \,.
\]
If this is the case, we write $\reallywidehat{\mu} \deq \nu$ and call it the \textit{Fourier transform of $\mu$ in the strongly tempered sense}.  

We write $\STM(\RR^d)$ for the space of all strongly tempered measures $\mu$ which are Fourier transformable in the strongly tempered sense.
\end{definition}

In the rest of the paper, all measures are assumed to be strongly tempered. Whenever we say that a measure $\mu$ is Fourier transformable, it is understood to mean Fourier transformable in the strongly tempered sense. Let us make some important remarks.

\begin{remark}\label{r1}
\begin{itemize}
\item[(a)]Let $\mu$ be a strongly tempered measure. Then $\mu$ is Fourier transformable in the measure sense (see \cite{AG1974,MS2017} for a definition) if and only if $\mu\in\STM(\RR^d)$ and $\reallywidehat{\mu}$ is translation bounded \cite{SS2021}. Moreover, in this case the two Fourier transforms coincide.

In particular, if $\mu$ is a strongly tempered measure that is also Fourier transformable in the measure sense, then $\mu\in\STM(\RR^d)$ and the Fourier transforms coincide. This allows us to use the same notation for both Fourier transforms. 
\item[(b)] Typically, the diffraction measure $\reallywidehat{\gamma}$ is defined as the Fourier transform in the measure sense of the autocorrelation measure $\gamma$. In general, $\gamma$ is a strongly tempered measure. Therefore, by (a), the diffraction measure $\reallywidehat{\gamma}$ is the Fourier transform in the strongly tempered sense of $\gamma$.
\end{itemize}
\end{remark}

We now list some easy observations about $\STM(\RR^d)$:

\begin{fact}
Let $\mu\in \STM(\RR^d)$. Then the following statements hold: 
\begin{itemize}
    \item[(a)] The measures $\mu$ and $\reallywidehat{\mu}$ uniquely determine each other.
    \item[(b)] For all $f \in \cS(\RR^d)$ we have $\hatf \in L^1(|\mu|)$, $f\in L^1(|\reallywidehat{\mu}|)$, and
    \[
    \int_{\RR^d} \reallywidehat{f}(y) \dd \mu(y) = \int_{\RR^d} f(x) \dd \reallywidehat{\mu}(x) \,.
    \]
    \item[(c)] $\reallywidehat{\mu} \in \STM(\RR^d)$ and 
    \[
    \reallywidehat{\reallywidehat{\mu}} =\mu^\dagger\,,
    \]
    where $\mu^\dagger$ is the measure defined by $\mu^\dagger(f):= \int_{\RR^d} f(-x)\, \dd \mu(x)$.

    In particular, the Fourier transform is a bijection from $\STM(\RR^d)$ onto $\STM(\RR^d)$. 
\end{itemize}
\end{fact}

The results in this paper rely on these observations.  
%In particular, the work is largely motivated by the fact (b), since it allows us to relate the measure of a circle to the autocorrelation via the Fourier transform of Schwartz functions (see Section~\ref{sec:circles}). When $\gamma$ is the autocorrelation of a translation bounded measure,  (b) also holds for any function in the Feichtinger algebra (see \cite{Jak2018}), so the results Section~\ref{sec:circles} can be generalized in this direction. To keep the level of technicality in the paper low,  we focus on Schwartz functions.  
 
Let us conclude this section by recalling a classical result about convolutions that we will need in Section~\ref{sec:diffraction}. 

\begin{theorem}[{\cite[Theorem 2.3.20]{Gra2008}}]\label{thm:Gra08}
    Let $T\in\cS'(\RR^d)$  be a tempered distribution, and let $f\in\cS(\RR^d)$ be a Schwartz function. Then $f * T$ is a $C^\infty$ function, where
    \[
    f * T(\varphi) := T(f^\dagger * \varphi) \qquad \forall \varphi \in \cS(\RR^d) \,,
    \]
    and $f^\dagger$ denotes the function given by $f^\dagger(x)=g(-x)$ for all $x\in\RR^d$.
    Furthermore, if $T$ has compact support, then $f * T$ is a Schwartz function.\qed
\end{theorem}

\begin{remark} Any measure of compact support is a strongly tempered measure. Moreover, 
     a strongly tempered measure $\mu$ has compact support if and only if the tempered distribution $T_\mu$ is compactly supported as a distribution.
\end{remark}

%\begin{fact}
%For all $f \in \cS(\RR^d)$ and all $\mu \in \STM(\RR^d)$ we have $\hatf \in L^1(|\mu|)$, $f\in L^1(|\hatmu|)$, and
%    \[
%    \int_{\RR^d} \reallywidehat{f}(y) \dd \mu(y) = \int_{\RR^d} f(x) \dd \, \reallywidehat{\mu} (x) \,.
%    \]
%\end{fact}

\section{Bessel functions}\label{sec:Bessel}

In Section~\ref{sec:diffraction} we are going to run into certain integrals. It turns out that these are exactly the integrals defining the so called Bessel function $J_0$ of first kind and order zero, and the modified Bessel function $I_0$ of the first kind and order zero. The asymptotic behaviour of $I_0$ also plays an important role in our computations. In this section we collect the results about these functions that will be important for us.

%In this section, we collect facts about Bessel functions that will be essential for the results in Section~\ref{sec:rcircle} and Section~\ref{sec:diffraction}.

\subsection{$J_0$ and the Fourier transform of circles}

\begin{definition}\label{def:J0}
    For each $t\geq 0$, define
\[
    J_0(t) \deq \frac{1}{2\pi}\int_{-\pi}^\pi e^{-it\cos\theta}\,\dd\theta\,.
\]
    $J_0$  is called the \emph{Bessel function of the first kind and order 0}.
\end{definition}\label{def:thetar}
Since $\cos\theta$ is $2\pi$-periodic, it is immediate that for all $a \in \RR$ we have
\[
\frac{1}{2\pi}\int_{a}^{a+\pi} e^{-it\cos\theta}\,\dd\theta = J_0(t) \,.
\]
\medskip

The Bessel function $J_0$ will be important for the following reason: the Fourier transform of the uniform probability measure of a circle is given by a $J_0$ function.  Let us make this statement precise.

\begin{definition}
    For $r>0$, define $\theta_r:\Cc(\RR^2)\rightarrow \CC$ by 
    $$
    \theta_r(f) = \frac{1}{2\pi}\int_{-\pi}^\pi f(r\cos\theta, r\sin\theta)\,\dd\theta \qquad \forall f\in\Cc(\RR^2)\,.
    $$
    Then, $\theta_r$ is a probability measure with compact support $C_r:= \{x \in \RR^2 : \| x \|=r \}$.
\end{definition}

\begin{remark} Note that $C_r$ is a compact Abelian group under the group operation 
\[
(re^{i \theta}) \diamond (re^{i \varphi}) := re^{i (\theta+\varphi)}  \,.
\]
Then, $\theta_r$ can be seen as the unique probability Haar measure on $(C_r,\diamond)$, embedded into $\RR^2$.  
\end{remark}

Since $\theta_r$ has compact support, it defines a strongly tempered measure, and we can easily compute its Fourier transform using \eqref{eq:TD-FT}. 

\begin{proposition}\label{prop:circle-FT}
    Let $r>0$.  The Fourier transform of $\theta_r$ is an absolutely continuous measure with density function
    \[
    \reallywidehat{\theta_r}(y) =J_0(2\pi r \|y\|) \qquad \forall y\in\RR^2\,.
    \]
\end{proposition}
\begin{proof}
Since $\theta_r$ is a distribution with compact support, its Fourier transform is an absolutely continuous measure with density function given by
    \begin{align*}
\reallywidehat{\theta_r}(y) &= \int_{\RR^d} e^{-2\pi i x \cdot y}\,\dd \theta_r(x)  = \frac{1}{2\pi}\int_{-\pi}^\pi  e^{-2\pi i r (y_1 \cos\theta+y_2 \sin\theta) } \dd \theta \qquad \forall y=(y_1,y_2)\in\RR^2\,.
\end{align*}
Consider $y = (R\cos \phi, R\sin \phi)$ for some $R\geq 0$ and $\phi\in[-\pi,\pi)$.  We have
\begin{align*}
\reallywidehat{\theta_r}(R\cos(\phi), R \sin(\phi)) &= \frac{1}{2\pi}\int_{-\pi}^\pi  e^{-2\pi i rR \cos(\phi+\theta) } \dd \theta \\
&\stackrel{\psi=\phi+\theta}{=\joinrel=\joinrel=\joinrel=\joinrel=} \frac{1}{2\pi}\int_{\phi-\pi}^{\phi+\pi}  e^{-2\pi i rR \cos\psi } \dd \psi = J_0(2\pi rR)\,.
\end{align*}
This shows that $\reallywidehat{\theta_r}(y) = J_0(2\pi r\|y\|)$ for all $y\in\RR^2$, as desired.
\end{proof}

\subsection{$I_0$ and its asymptotic behaviour }

\begin{definition}\label{def:I0}
    For each $t\geq 0$, define
    $$
    I_0(t) \deq \frac{1}{2\pi}\int_{-\pi}^\pi e^{t\cos\theta}\,\dd\theta\,.
    $$
    $I_0$ is called the \emph{modified Bessel function of the first kind of order 0}.
\end{definition}

 The function $I_0$ will appear naturally in Section~\ref{sec:diffraction}  as the convolution of a Gaussian function and $\theta_r$.  This type of convolution will be essential to our approach. The following asymptotic property of $I_0$ will be crucial in our computations:

\begin{theorem}[Asymptotic behavior of $I_0$]\label{thm:I0-asymptotic}
\[
\lim_{t\rightarrow\infty} \frac{\sqrt{2\pi t}I_0(t)}{e^{t}} =  1\,.
\]
\end{theorem}

This result is well-known and follows from a standard application of Laplace's method for asymptotic expansions of functions of the form 
\[
I(x) = \int_a^b e^{xg(t)}f(t)\,\dd t\,,
\]
(refer to \cite[Ch.\,6.4]{BO1999} for the general theory of Laplace's method).
The limit can easily be deduced from this general theory; 
%{\color{green} however the proofs are typically done for the entire family of modified Bessel functions of first kind, and are technical} \red{(actually, the proofs give the asymptotic expansion of $I(x) = \int_a^b e^{x g(t)} f(t)$ in terms of the asymptotic expansion of $f(t)$ as $t\rightarrow 0^+$, and the $I_n$ functions are usually written as an example, where the general theorem is applied.) }.  Thus, 
anyhow, for completeness, we include a short direct proof of Theorem~\ref{thm:I0-asymptotic} below.

First, we need the following lemma.

\begin{lemma}\label{lem:I0-integral}
    For each $t\geq 0$, we have
    $$
    I_0(t) = \frac{e^t}{\pi}\int_0^2 \frac{e^{-t\xi}}{\sqrt{\xi(2-\xi)}}\,\dd\xi\,.
    $$
\end{lemma}

\begin{proof}
    Using symmetry and applying the substitution $\xi = 1-\cos\theta$, we have
    \begin{align*}
    I_0(t) &= \frac{1}{2\pi}\int_{-\pi}^\pi e^{t\cos\theta}\,\dd\theta =  \frac{1}{\pi}\int_0^\pi e^{t\cos\theta}\,\dd\theta 
    \stackrel{\xi = 1-\cos\theta}{=\joinrel=\joinrel=\joinrel=\joinrel=\joinrel=\joinrel=} \frac{1}{\pi}\int_0^2 \frac{e^{t(1-\xi)}}{\sqrt{\xi(2-\xi)}}\,\dd\xi \\
    &= \frac{e^t}{\pi}\int_0^2 \frac{e^{-t\xi}}{\sqrt{\xi(2-\xi)}}\,\dd\xi\,.
    \end{align*}
\end{proof}

We can now prove the asymptotic result.

\begin{proof}[Proof of Theorem~\ref{thm:I0-asymptotic}]
By Lemma~\ref{lem:I0-integral}, we have
\begin{align*}
\frac{\pi I_0(t)}{e^{t}}&=  \int_0^2  \frac{e^{-t\xi}}{\sqrt{\xi(2-\xi)}}\,\dd \xi \\
&=\underbrace{\int_0^2  \frac{e^{-t\xi}}{\sqrt{2\xi}}\,\dd \xi}_{I_1}+\underbrace{\int_0^1  \frac{e^{-t\xi}}{\sqrt{\xi}} \left( \frac{1}{\sqrt{2-\xi}} - \frac{1}{\sqrt{2}}\right)\,\dd \xi}_{I_2}+\underbrace{\int_1^2  \frac{e^{-t\xi}}{\sqrt{\xi}} \left( \frac{1}{\sqrt{2-\xi}} - \frac{1}{\sqrt{2}}\right)\,\dd \xi}_{I_3} \,.
\end{align*}
We now look at each of the three integrals one by one.  First, with a simple substitution we get
\begin{align*}
I_1&\stackrel{u=\sqrt{t \xi}}{=\joinrel=\joinrel=\joinrel=\joinrel=}  \sqrt{\frac{2}{t}} \int_0^{\sqrt{2t}}  e^{-u^2}
\dd u \,.
\end{align*}
Then, using $\int_{-\infty}^\infty e^{-u^2}\,\dd u = \sqrt{\pi}$, we get
\begin{align*}
\left|  \sqrt{\frac{2t}{\pi}} I_1-1 \right| = \left|  \frac{2}{\sqrt{\pi}} \int_0^{\sqrt{2t}}   e^{-u^2}
\dd u-1 \right| &= \frac{2}{\sqrt{\pi}} \int_{ \sqrt{2t}}^{\infty} e^{-u^2} \dd u \\
&\leq \frac{2}{\sqrt{\pi}} \int_{ \sqrt{2t}}^{\infty} e^{-u\sqrt{2t}} \dd u =
\frac{2}{\sqrt{\pi}}\frac{e^{-2t}}{\sqrt{2t}}\,.  
\end{align*}
Next, observe that the function $f: [0,1] \to \RR $
\[
f(\xi):=\left\{
\begin{array}{cc}
\frac{1}{\sqrt{\xi}}\lt( \frac{1}{\sqrt{2-\xi}} - \frac{1}{\sqrt{2}}\rt)     & \mbox{ if } \xi \neq 0   \\
    0 & \mbox{ if } \xi =0\,,
\end{array}
\right.
\]
is continuous and hence bounded on $[0,1]$. Therefore, there exists some $M>0$ such that
\begin{align*}
 |I_2|& \leq \int_0^1 M e^{-t\xi}\,\dd \xi = \frac{M}{t} (1-e^{-t})\,.
\end{align*}
For the third integral, we use that $e^{-t\xi}\leq e^{-t}$ for $t\geq 1$ to obtain 
\begin{align*}
|I_3|& \leq e^{-t} \int_1^2  \frac{1}{\sqrt{\xi}} \left| \frac{1}{\sqrt{2-\xi}} - \frac{1}{\sqrt{2}}\right|\,\dd \xi = Ce^{-t}\,,
\end{align*}
where 
\[
C:= \int_1^2  \frac{1}{\sqrt{\xi}} \left| \frac{1}{\sqrt{2-\xi}} - \frac{1}{\sqrt{2}}\right|\,\dd \xi < \infty\,.
\]
Putting everything together, we get 
\begin{align*}
\left|\frac{\sqrt{2\pi t}I_0(t)}{e^{t}} - 1  \right| &=\left|  \sqrt{\frac{2t}{\pi}}( I_1+I_2+I_3)-1 \right| \\
&\leq \left|  \sqrt{\frac{2t}{\pi}} I_1-1 \right| + \sqrt{\frac{2t}{\pi}}\left( |I_2|+|I_3| \right) \\
&\leq \frac{2}{\sqrt{\pi}}\frac{e^{-2t}}{\sqrt{2t}} + \sqrt{\frac{2t}{\pi}}\lt(\frac{M}{t} (1-e^{-t})+Ce^{-t}\rt)  \,.
\end{align*}
The claim follows.
\end{proof}

\section{The measure of circles}\label{sec:circles}

We now restrict our attention to $\RR^2$. In this section, we show that the measure of a circle can computed as a limit via a sequence of Schwartz functions. For a measure $\mu$ which is Fourier transformable, this will allow us to express the measure $\reallywidehat{\mu}(C_r)$ as a formula in terms of the original measure $\mu$. 

We will more generally show that we can use any sequence of Schwartz functions that satisfy a condition, which we will call and an approximate $r$-circle. While approximate $r$-circles are easy to construct, for applications we will need such a sequence with nice, easily computable Fourier transforms. We will construct one such example in Section~\ref{sec:diffraction}.

We now introduce the notion of an approximate $r$-circle.  

\begin{definition}\label{def:r-circle}
    Let $r \geq 0$. We say that a sequence $\seq{f_n}$ in $\cS(\RR^2)$ is an \textit{approximate $r$-circle} if the following conditions are satisfied:
    \begin{enumerate}[label=(C\arabic*)]
        \item \label{itm:C1} $\lim\limits_{n\rightarrow\infty}f_n(x) = 1_{C_r}(x)$ for all $x\in\RR^2$;
        \item \label{itm:C2} there exists some $F\in \cS(\RR^2)$ such that $$|f_n(x)|\leq F(x)\qquad \forall x\in\RR^2\,,\forall n\in\NN\,.$$
    \end{enumerate}
\end{definition}

In the following lemma, we see that for any measure in $\STM(\RR^2)$, the measure of $C_r$ can be computed as a limit via an approximate $r$-circle. 

\begin{lemma}\label{lem:DCT}
    Let $\mu$ be a strongly tempered measure. 
 Then for any approximate $r$-circle $\seq{f_n}$, we have
     \begin{equation}\label{eq:DCT}
     \mu(C_r) = \lim_{n\rightarrow\infty} \int_{\RR^2} f_n(x)\,\dd \mu(x) \,.
     \end{equation}
\end{lemma}

\begin{proof}
    By Fact~\ref{fact:strongly-tempered}, if $\mu$ is strongly tempered then $f\in L^1(|\mu|)$ for every $f\in\cS(\RR^2)$.  In particular, if $\seq{f_n}$ is an approximate $r$-circle and $F$ is as in condition \ref{itm:C2}, then $F\in L^1(|\mu|)$.  Then by condition \ref{itm:C1} and the dominated convergence theorem, we have
    $$
    \mu(C_r) = \int_{\RR^2} 1_{C_r}(x)\,\dd\mu(x) = \int_{\RR^2} \lim_{n\rightarrow\infty} f_n(x) \,\dd\mu(x) = \lim_{n\rightarrow\infty} \int_{\RR^2} f_n(x)\,\dd\mu(x)\,.
    $$
\end{proof} 

\begin{remark}
    Using the monotone convergence theorem, we also get that \eqref{eq:DCT} holds for any decreasing sequence $\seq{f_n}$ of functions in $L^1(|\mu|)$ satisfying \ref{itm:C1}.
\end{remark}

As an immediate consequence of Lemma~\ref{lem:DCT} we get the following result, which allows us calculate the intensity of a circle in the diffraction in terms of the autocorrelation measure. 

\begin{corollary}\label{cor:STM}
    If $\mu\in\STM(\RR^d)$ and $\seq{f_n}$ is an approximate $r$-circle, then 
    \[
    \hatmu(C_r) = \lim_{n\rightarrow\infty} \int_{\RR^2} f_n(x)\,\dd\,\hatmu(x) = \lim_{n\rightarrow\infty} \int_{\RR^2} \reallywidehat{f_n}(x)\,\dd\mu(x)\,. 
    \qed\]
\end{corollary}

This formula will become useful once we construct examples of approximate $r$-circles with functions that have nice Fourier transforms, which we do in the next section.

\section{Circles in diffraction}\label{sec:diffraction}
%\section{An approximate $r$-circle} \label{sec:rcircle}

In this section, we construct an example of an approximate $r$-circle via the convolution of $\theta_r$ with a family of Gaussian functions.  We then use this approximate $r$-circle and Corollary~\ref{cor:STM} to derive an explicit formula for the diffraction measure of a circle in terms of the autocorrelation.  A different example can be found in Section~\ref{sect:anotherapp}.

Let $G_n$ be the following family of two dimensional Gaussian functions 
\begin{equation}\label{eq:alt-Gn}
G_n(x)=e^{-n^2\|x\|^2} \qquad \forall x\in\RR^2\,, \ \forall n\in\NN\,, 
\end{equation}
and for each $r>0$, let $\theta_r$ be the uniform probability measure on the circle $C_r$.  Consider the convolution $G_n*\theta_r$, defined by
$$
G_n*\theta_r(x) = \frac{1}{2\pi}\int_{-\pi}^\pi G_n(x_1 - r\cos\theta, x_2 - r\sin\theta)\,\dd\theta \qquad \forall x = (x_1,x_2)\in\RR^2\,.
$$
Since both $\theta_r$ and $G_n$ are rotationally invariant, their convolution $G_n*\theta_r$ is also rotationally invariant. 
 With this observation, it is easy to show the following.
 
 \begin{lemma}\label{lem:alt-conv}
     For each $r>0$ and $n\in\NN$, we have
     \begin{equation}\label{eq:alt-conv}
        G_n*\theta_r(x) =  I_0(2n^2r\|x\|)e^{-n^2(\|x\|^2 + r^2)} \qquad \forall x\in\RR^2\,.
    \end{equation}
 \end{lemma}

 \begin{proof}
 Using rotational invariance of $G_n*\theta_r$, for each $x\in\RR^2$ we have:
 \begin{align*}
    & G_n*\theta_r(x) = G_n*\theta_r(\|x\|,0) = \frac{1}{2\pi}\int_{-\pi}^\pi e^{-n^2((\|x\|-r\cos\theta)^2 + r^2\sin^2\theta)}\,\dd\theta \\
    &= \frac{1}{2\pi} e^{-n^2(\|x\|^2 + r^2)} \int_{-\pi}^\pi e^{2n^2r\|x\|\cos\theta} \,\dd\theta 
    %= \frac{1}{\pi} e^{-n^2(\|x\|^2 + r^2)} \int_0^\pi e^{2n^2r\|x\|\cos\theta} \,\dd\theta \\
    %&= e^{-n^2(\|x\|^2 + r^2)} I_0(2n^2r\|x\|) 
    = I_0(2n^2r\|x\|)e^{-n^2(\|x\|^2 + r^2)}\,.
 \end{align*}
 \end{proof}

 Throughout the rest of this section, we fix $r>0$ and define the functions
$$
g_n \deq 2nr\sqrt{\pi} G_n \,, \qquad f_n = g_n*\theta_r = 2nr\sqrt{\pi} G_n*\theta_r\,.
$$
We will prove that $f_n$ is an approximate $r$-circle. Note that by \eqref{eq:alt-conv} we have 
\[
f_n(x)= 2nr\sqrt{\pi} I_0(2n^2r\|x\|)e^{-n^2(\|x\|^2 + r^2)} \qquad \forall x\in\RR^2\,.
\]

First, we show that $f_n$ converges pointwise to $1_{C_r}$.

\begin{lemma}\label{lem:rcircle}
 For each $r>0$, the sequence $\seq{f_n}$ is an approximate $r$-circle. 
\end{lemma}
\begin{proof}
Let us first note that $f_n\in\cS(\RR^2)$ for each $n\in\NN$ by Theorem~\ref{thm:Gra08}, since $G_n$ is a Schwartz function and $\theta_r$ is compactly supported.

Next, observe that
\begin{equation} \label{eqx}
\begin{split}
f_n(x) &=\frac{2nr\sqrt{\pi}I_0(2n^2 r\|x\|)}{e^{n^2(\|x\|^2 + r^2)}}=\frac{2n\sqrt{r\pi \| x\|}I_0(2n^2 r\|x\|)}{e^{2 n^2 r \|x\|}} \cdot  \frac{\sqrt{r} e^{2 n^2 r \|x\|}} {\sqrt{\|x\|}e^{n^2(\|x\|^2 + r^2)}} \\
&= u(2n^2r\|x\|) \sqrt{\frac{r}{\|x\|}} e^{-n^2(\|x\|-r)^2}
\end{split}
\end{equation}
where 
\[
u(t)\deq \frac{\sqrt{2\pi t}I_0(t)}{e^{t}} \,.
\]
We now can show that the two conditions (C1) and (C2) of Definition~\ref{def:r-circle} hold.

\textbf{(C1):}
Since $I_0(0)=1$, when $x = 0$ we have
\[
\lim_{n\rightarrow\infty} f_n(0) = 2nr\sqrt{\pi}e^{-n^2r^2} = 0\,.
\]
Next, by Theorem~\ref{thm:I0-asymptotic} and \eqref{eqx}, for all $x \neq 0$ we have 
\[
\lim_{n\rightarrow\infty} f_n(x)=\lim_{n\rightarrow\infty}\sqrt{\frac{r}{\|x\|}} e^{-n^2(\|x\|-r)^2}= 
    \begin{cases}
        1 & \|x\| = r \\
        0 & \|x\|\neq 0 \ \text{and} \ \|x\| \neq r\,.
    \end{cases} \,.
\]
This proves (C1).

\textbf{(C2):} The function $u$ is continuous on $[0, \infty)$ and satisfies 
\[
\lim_{t\rightarrow\infty} u(t) = 1\,.
\]
This implies that there exists some constant $A$ such that $|u(t)| \leq A$ for all $t \in [0, \infty)$.
In particular, 
\[
|u(2n^2r\|x\|)| \leq A \qquad \forall n\in\NN, \forall x \in \RR^2 \,. 
\]
By \eqref{eqx} we then have 
\[
|f_n(x)| \leq A \sqrt{\frac{r}{\|x\|}} e^{-n^2(\|x\|-r)^2} \leq A \sqrt{\frac{r}{\|x\|}} e^{-(\|x\|-r)^2} =: \psi_1(x) \,.
\]
In particular, we get 
\begin{equation}\label{eqx1}
|f_n(x)| \leq A \sqrt{2} e^{-(\|x\|-r)^2} \qquad \forall \| x \| \geq \frac{r}{2} \,.   
\end{equation}

\noindent Next, let $x \in B_{\frac{r}{2}}$.  Then, for any $y\in C_r$ we have 
        \[
        \|x-y\|\geq r - \|x\| \geq \frac{r}{2}\,.
        \] 
        From this, we see that
        \begin{align*}
        0 \leq g_n*\theta_r(x) &= \int_{\RR^2} g_n(x-y) \,\dd\theta_r(y) = \int_{C_r} g_n(x-y) \,\dd\theta_r(y) \\
        &= 2rn\sqrt{\pi}\int_{C_r} e^{-n^2\|x-y\|^2}\,\dd \theta_r(y) \\
        &\leq 2rn\sqrt{\pi}\int_{C_r} e^{-n^2r^2/4}\,\dd \theta_r(y) \\
        &= 2rn\sqrt{\pi}e^{-n^2r^2/4} \leq 2r\sqrt{\pi}\cdot c_1 =: c_2\,,
        \end{align*}
        where
        \[
        c_1 \deq \sup \{ ne^{-n^2r^2/4} : n \in \NN \} < \infty\,.
        \]
Now, pick some $\psi_2 \in \Cc^\infty(\RR^2)$ such that 
\[
\psi_2 \geq c_2 1_{B_{\frac{r}{2}}} \,.
\]
Then, we have 
\begin{equation}\label{eqx2}
|f_n(x)| \leq \psi_2(x)  \qquad \forall \| x \| \leq \frac{r}{2} \,.   
\end{equation}
Therefore, combining \eqref{eqx1} and \eqref{eqx2} we get 
\[
|f_n(x)| \leq \psi_1(x)+\psi_2(x)  \qquad \forall \in \RR^2 \,.
\]
Setting $F:=\psi_1+\psi_2 \in \cS(\RR^2)$, we showed condition (C2).
\end{proof}

We have thus constructed an explicit approximate $r$-circle.  Since we can easily compute $\reallywidehat{f_n}$, we can now construct an explicit formula for the diffraction measure $\hatgam(C_r)$ for a circle $C_r$ of a given radius $r>0$, in terms of the autocorrelation $\gamma$.  To derive the formula, we simply need to apply Corollary~\ref{cor:STM} for the approximate $r$-circle $\seq{f_n}$, and compute $\reallywidehat{f_n}$.  

We proved in Lemma~\ref{lem:rcircle} that $\seq{f_n}$ is an approximate $r$-circle, so by Corollary~\ref{cor:STM}, we have the formula
\[
\hatmu(C_r) = \lim_{n\rightarrow\infty} \int_{\RR^2} f_n(x) \,\dd \, \hatmu(x) = \lim_{n\rightarrow\infty} \int_{\RR^2} \reallywidehat{f_n}(y)\,\dd\mu(y)\,,
\]
for any measure $\mu\in\STM(\RR^2)$. It is easy to see that
\[
\reallywidehat{G_n}(y) = \frac{\pi}{n^2}e^{-\pi^2\|y\|^2/n^2} \qquad \forall y\in\RR^2\,.
\]
Moreover, by Proposition~\ref{prop:circle-FT} we have
\[
\reallywidehat{\theta_r}(y)=J_0(2\pi r \|y\|) \qquad \forall y\in\RR^2\,.
\]
Therefore, by the convolution theorem we get:
\begin{lemma} The function 
\[
f_n = 2nr\sqrt{\pi}G_n*\theta_r= 2nr\sqrt{\pi} I_0(2n^2r\|x\|)e^{-n^2(\|x\|^2 + r^2)} \,,
\]
is an approximate $r$-circle with Fourier transform given by
\[
\reallywidehat{f_n}(y) = \frac{2r\pi\sqrt{\pi}}{n}e^{-\pi^2\|y\|^2/n^2}J_0(2\pi r \|y\|) \qquad \forall y \in \RR^2\,.
\]
\qed
\end{lemma}

Combining all results from this paper we get:

\begin{theorem}\label{thm:diffraction-formula}
    Let $\mu \in \STM(\RR^2)$ and let $r>0$. Then
    \begin{equation}\label{eq:diffraction-formula}
        \hatmu(C_r) = \lim_{n\rightarrow\infty} \frac{2r\pi\sqrt{\pi}}{n}\int_{\RR^2} e^{-\pi^2\|y\|^2/n^2}J_0(2\pi r \|y\|) \,\dd\mu(y)\,.
    \end{equation} \qed
\end{theorem}

\begin{remark}
   When $\mu=\gamma$ is an autocorrelation measure, the above formula calculates the diffraction measure $\hatgam(C_r)$ for the circle $C_r$ in terms of the autocorrelation measure $\gamma$.
\end{remark}

\section{Applications}\label{sec:applications}

In this section, we provide applications of Theorem~\ref{thm:diffraction-formula} and highlight some interesting observations that could lead to further applications.  Below, we use $\lambda$ to denote the Lebesgue measure on $\RR^2$, and sometimes use the notation $\dd\lambda(x)$ instead of $\dd x$ in cases where the integral could be mistaken for a line integral.

We start with some simple examples with measures $\mu$ where $\hatmu(C_r)$ is known and easily evaluated. 
%\begin{example}
%    Let $\mu=\delta_0$ be the Dirac measure concentrated at the origin.  Then $\hatmu = \lambda$ is the Lebesgue measure on $\RR^2$, and \eqref{eq:diffraction-formula} becomes
%    $$
%    \lim_{n\rightarrow\infty} \frac{2r\pi\sqrt{\pi}}{n}\int_{\RR^2} e^{-\pi^2\|y\|^2/n^2}J_0(2\pi r \|y\|) \,\dd \delta_0(y) = \lim_{n\rightarrow\infty} \frac{2r\pi\sqrt{\pi}}{n} = 0\,.
%    $$
%\end{example}

\begin{example}\label{ex:deltax}
Let $\mu = \delta_x$ for some $x\in\RR^2$.  Then $\hatmu(y) = e^{-2\pi i x\cdot y}$ and \eqref{eq:diffraction-formula} becomes
\[
\int_{C_r} e^{-2\pi i x\cdot y} \,\dd \lambda(y) = \lim_{n\rightarrow\infty} \frac{2r\pi\sqrt{\pi}}{n} e^{-\pi^2 \|x\|^2/n^2} J_0(2\pi r\|x\|) = 0\,,
\]
where $\lambda$ is the Lebesgue measure in $\RR^2$.
\end{example}

\begin{example}
Let $\mu = \theta_{r'}$.  Recall from Proposition~\ref{prop:circle-FT} that for any $r'>0$, we have that
\[
\reallywidehat{\theta_{r'}}(y) = J_0(2\pi r' \|y\|) \qquad \forall y \in \RR^2\,.
\]
Thus, in this case, \eqref{eq:diffraction-formula} gives us:
\[
\int_{C_r} J_0(2\pi r'\|x\|)\,\dd \lambda(x) 
= \lim_{n\rightarrow\infty} \frac{r\pi^{\frac{1}{2}}}{n}\int_{-\pi}^\pi e^{-\pi^2(r')^2/n^2}J_0(2\pi r r') \,\dd \theta = 0\,.
\]
\end{example}

\begin{remark} More generally, if $\mu \in \STM(\RR^2)$ is such that $\reallywidehat{\mu}$ is an absolutely continuous measure, then 
\[
 \lim_{n\rightarrow\infty} \frac{2r\pi\sqrt{\pi}}{n}\int_{\RR^2} e^{-\pi^2\|y\|^2/n^2}J_0(2\pi r \|y\|) \,\dd\mu(y) =0 \,.
\]
\end{remark}

Next, let us look at the cases where $\mu$ is a character.
\begin{example}
    Let $\mu = \lambda$ be the Lebesgue measure on $\RR^2$.  Then $\hatmu = \delta_0$ and \eqref{eq:diffraction-formula} tells us that for all $r>0$ we have 
    \[
\lim_{n\rightarrow\infty} \frac{2r\pi\sqrt{\pi}}{n}\int_{\RR^2} e^{-\pi^2\|y\|^2/n^2}J_0(2\pi r \|y\|) \,\dd y = 0\,.
    \]
\end{example}

\begin{example}
    Let $x\in \RR^2\backslash \{0\}$ and let $\mu$ be the function $y \mapsto e^{-2\pi i x\cdot y}$.  Then $\hatmu = \delta_x$ and \eqref{eq:diffraction-formula} becomes
    $$
 \lim_{n\rightarrow\infty} \frac{2r\pi\sqrt{\pi}}{n}\int_{\RR^2} e^{-2\pi i x\cdot y}e^{-\pi^2\|y\|^2/n^2}J_0(2\pi r \|y\|) \,\dd y 
    = \begin{cases}
        1 & \|x\| = r \\
        0 & \|x\|\neq r\,.
    \end{cases}
    $$
    Notice that the above limit is just $\lim_{n\rightarrow\infty}f_n(x)$, where $f_n$ is the approximate $r$-circle defined in Section~\ref{sec:diffraction}. 
\end{example}

Next, we consider the case where $\hatmu = \theta_{r'}$ for some $r'>0$.  This example is of particular interest because it leads to an orthogonality relation for Bessel functions of the form $J_0(2\pi r\|\cdot\|)$.
Applying Corollary~\ref{cor:STM} with $\mu = J_0(2\pi r\|\cdot\|)$, we get the following:

\begin{theorem}
    Let $r, r'>0$.  If $\seq{f_n}$ is any approximate $r$-circle, then the following identity holds:
    \[ \lim_{n\rightarrow\infty}\frac{2r\pi\sqrt{\pi}}{n}\int_{\RR^2} \reallywidehat{f_n}(y) J_0(2\pi r'\|y\|) \,\dd y = \begin{cases}
1 & r' = r \\
0 & r' \neq r\,.
\end{cases}
    \]
    In particular, when $f_n = 2nr\sqrt{\pi}G_n*\theta_r$ we have
    \begin{align}
        &\lim_{n\rightarrow\infty} \frac{r\sqrt{\pi}}{n}\int_{0}^\infty e^{-z^2/4n^2}J_0(r z) J_0(r'z) z \,\dd z  \label{eq:orthogonal}\\
        &= \lim_{n\rightarrow\infty} \frac{4r\pi^2\sqrt{\pi}}{n}\int_{0}^\infty e^{-\pi^2 u^2/n^2}J_0(2\pi r u) J_0(2\pi r'u) u \,\dd u \nonumber \\
        &=\lim_{n\rightarrow\infty} \frac{2r\pi\sqrt{\pi}}{n}\int_{\RR^2} e^{-\pi^2\|y\|^2/n^2}J_0(2\pi r \|y\|) J_0(2\pi r'\|y\|) \,\dd y= 
        \begin{cases}
        1 & r' = r \\
        0 & r' \neq r\,.
        \end{cases}\nonumber
    \end{align}
\end{theorem}

\begin{proof}
Corollary~\ref{cor:STM} applied to $\mu = J_0(2\pi r\|\cdot\|)$ gives
\[
\lim_{n\rightarrow\infty} \frac{2r\pi\sqrt{\pi}}{n}\int_{\RR^2} e^{-\pi^2\|y\|^2/n^2}J_0(2\pi r \|y\|) J_0(2\pi r'\|y\|) \,\dd y= 
        \begin{cases}
        1 & r' = r \\
        0 & r' \neq r\,.
        \end{cases} \,.
\]
Next, switching to polar coordinates gives
\[
\int_{\RR^2} e^{-\pi^2\|y\|^2/n^2}J_0(2\pi r \|y\|) J_0(2\pi r'\|y\|) \,\dd y= 2 \pi \int_{0}^\infty e^{-\pi^2 u^2/n^2}J_0(2\pi r u) J_0(2\pi r'u) u \,\dd u \,.
\]
The substitution $z=2 \pi u$ yields the remaining formula.
\end{proof}

\begin{remark}
    The integrals that appear in \eqref{eq:orthogonal} can be evaluated using standard tables of integrals (e.g. \cite[6.633 (2)]{GR2000}). 
\end{remark}

The above examples show that Corollary~\ref{cor:STM} and Theorem~\ref{thm:diffraction-formula} can lead to interesting integral identities. In the next section, we will focus on pure point measures; however, before doing so, let us state here the following result, which follows immediately from the rotational invariance of the pinwheel diffraction \cite{MPS2006}.

\begin{proposition} Let $\gamma$ be the autocorrelation of the pinwheel tiling and let $r>0$. Then, there exists a circle of radius $r$ in the pinwheel diffraction if and only if 
\[
\lim_{n\rightarrow\infty} \frac{2r\pi\sqrt{\pi}}{n}\int_{\RR^2} e^{-\pi^2\|y\|^2/n^2}J_0(2\pi r \|y\|) \,\dd\gamma(y) \neq 0 \,.
\]\qed
\end{proposition} 

\begin{remark}
    The same also holds for the autocorrelation and diffraction of any tiling with statistical circular symmetry \cite{Fre2008}, as any such tiling has rotationally invariant diffraction.
\end{remark}

\subsection{Measures with pure point Fourier transform.}

Next we consider some examples with pure point diffraction.  In this case, the left side of \eqref{eq:diffraction-formula} is the sum of the intensities of the Bragg peaks that lie on the circle $C_r$.  More precisely, we get the following:

\begin{proposition}\label{prop:pp}
    Let $r>0$ and let $\mu\in\STM(\RR^2)$ be such that $\hatmu$ is pure point, i.e. $\hatmu$ is concentrated on a countable subset of $\RR^2$. Let $f_n$ be any approximate $r$-circle.
    Then
\[
\lim_{n\rightarrow\infty} \int_{\RR^2} \reallywidehat{f_n}(y) \,\dd\mu(y) =  \sum_{\substack{x\in \RR^2 \\ \| x \|=r} }\hatmu(\{x\})\,.
 \]

    In particular, when $f_n = 2nr\sqrt{\pi}G_n*\theta_r$ we have
\[
\lim_{n\rightarrow\infty} \frac{2r\pi\sqrt{\pi}}{n}\int_{\RR^2} e^{-\pi^2\|y\|^2/n^2}J_0(2\pi r \|y\|)\,\dd\mu(y) =  \sum_{\substack{x\in \RR^2 \\ \| x \|=r} }\hatmu(\{x\})\,.
 \] \qed
\end{proposition}

This result leads to an interesting formula for the sum of two squares function from number theory. Let us first recall its definition. 

\begin{definition} The \textbf{sum of two squares function} $r_2$ is the function $r_2: \NN \to \NN$ defined as 
\[
r_2(n):= \card \{ (k,l) \in \ZZ^2 : k^2+l^2= n \} \,.
\]
\end{definition}

\begin{remark} The sum of two squares function has the following properties:
\begin{itemize}
    \item[(a)] $r_2$ is an arithmetic function, meaning that $r_2(mn)=r_2(m)r_2(n)$ for all relatively prime $m,n$. 
    \item[(b)] $r_2$ has the alternate formula
    \[
    r_2(n)= 4(d_1(n)-d_3(n))\,,
    \]
    where $d_j(n)$ is the number of divisors of $n$ which equal $j \pmod{4}$.
    \item[(c)] In some places $r_2(n)$ is called the \textit{shelling function}, see for example  \cite{BGJR2000}, \cite{BG2003}.
\end{itemize}
\end{remark}

By the Poisson summation formula, we have that $\reallywidehat{\delta_{\ZZ^2}}=\delta_{\ZZ^2}$.  
Therefore, we get the following formula for $r_2(n)$ in terms of Bessel functions. We only list the result using the approximate $r$-circle from Section~\ref{sec:diffraction}, but similar formulas can be deduced for any approximate $r$-circle.

\begin{proposition}\label{prop:shelling} 
The sum of two squares function satisfies the following:
   \begin{itemize}
       \item[(a)] For all $k \in \NN$ we have 
       \[
r_2(k)= \lim_{n\rightarrow\infty} \frac{2\pi \sqrt{\pi k}}{n} \sum_{m =0}^{\infty}r_2(m) e^{-\frac{\pi^2m}{n^2}}J_0(2\pi \sqrt{km})  \,.
       \]
       \item[(b)] For all $r>0$ such that $r^2 \notin \NN$ we have 
       \[
 \lim_{n\rightarrow\infty} \frac{2r \pi\sqrt{\pi}}{n} \sum_{m =0}^{\infty}r_2(m) e^{-\frac{\pi^2m\textbf{}}{n^2}}J_0(2\pi  r \sqrt{m})  =0 \,.
       \]
   \end{itemize}
\end{proposition}
\begin{proof}
By Proposition~\ref{prop:pp}, for each $r>0$ we have
\[
\delta_{\ZZ^2}(C_r)= \lim_{n\rightarrow\infty} \frac{2r\pi\sqrt{\pi}}{n} \sum_{y \in \ZZ^2} e^{-\pi^2\|y\|^2/n^2}J_0(2\pi r \|y\|)\, \,.
\]
Next, grouping the elements in $\ZZ^2$ by their norm, we get
\begin{align*}
&\sum_{y \in \ZZ^2} e^{-\pi^2\|y\|^2/n^2}J_0(2\pi r \|y\|)= \sum_{m =0}^{\infty} \sum_{\substack{y\in\ZZ^2 \\ \|y\|^2=m}} e^{-\pi^2\|y\|^2/n^2}J_0(2\pi r \|y\|) \\
&= \sum_{m =0}^{\infty} \sum_{\substack{y\in\ZZ^2 \\ \|y\|^2=m}} e^{-\pi^2m/n^2}J_0(2\pi r \sqrt{m}) = \sum_{m =0}^{\infty}r_2(m) e^{-\pi^2m/n^2}J_0(2\pi r \sqrt{m}) \,.
\end{align*}
Now, if $r^2= k \in \NN$ we have 
\[
\delta_{\ZZ^2}(C_r)= r_2(k) \,,
\]
while when $r^2 \notin \NN$ we have 
\[
\delta_{\ZZ^2}(C_r)= 0 \,.
\]
\end{proof}

This can easily be extended to arbitrary lattices in $\RR^2$. Indeed, for a locally finite point set $\Lambda \subseteq \RR^2$ denote by $r_\Lambda : [0, \infty) \to \RR$ the function
\[
r_\Lambda(k):= \card \{ (s,t) \in \Lambda : s^2+t^2 =k \} =\delta_{\Lambda}(C_{\sqrt{k}}) \,. 
\]
Then, exactly as in Proposition~\ref{prop:shelling}  we have:
\begin{proposition} Let $L \subseteq \RR^2$ be a lattice with dual lattice $L^0$. Then, for all $k>0$ we have 
\[
r_L(k)=\frac{1}{\det(L)} \lim_{n\rightarrow\infty} \frac{2\pi\sqrt{\pi k}}{n} \sum_{m \in [0,\infty) }r_{L^0}(m) e^{-\frac{\pi^2m}{n^2}}J_0(2\pi  \sqrt{mk})\,.
\]\qed
\end{proposition}

\appendix
\section{A family of approximate $r$-circles}\label{sect:anotherapp}

In this appendix we offer an alternate approach to constructing approximate $r$-circles via convolution with the annulus 
%an annulus instead of a circle.  In particular, instead of using $\theta_r$, we use characteristic functions $1_{C_{r,\frac{1}{n}}}$, where
\[
C_{r,\frac{1}{n}}\deq \lt\{x\in\RR^2\,:\,r-\frac{1}{n} < \|x\| < r+\frac{1}{n}\rt\} \qquad \forall n\in\NN\,,
\]
instead of with $\theta_r$.

We provide general conditions on a sequence of functions $g_n$ that suffice for the convolution $g_n*1_{C_{r,\frac{1}{n}}}$ to be an approximate $r$-circle. The resulting formula for $\hatmu(C_r)$ is more complicated than that given in Theorem~\ref{thm:diffraction-formula}; however, the approach yields a large family of approximate $r$-circles, whereas the approach taken in Section~\ref{sec:diffraction} is specific to the convolution of $\theta_r$ with the Gaussian functions $G_n$.

\begin{lemma}\label{lem:alt-rcircle}
		Let $r > 0$, and let $\seq{g_n}$ be any sequence in $\cS(\RR^2)$ with the following properties:
		\begin{enumerate}[label=(\roman*)]
			\item $g_n \geq 0$ for all $n\in\NN$;
			\item $\lim\limits_{n\rightarrow\infty}\int_{B_{\frac{1}{n}}} g_n(x)\, \dd x = 1 $; 
			\item there exists some $\varphi\in\cS(\RR^2)$ such that 
			\[
			g_n(x)\leq \varphi(x) \qquad \forall \|x \| \geq \frac{1}{n} \,.
   \] 
		\end{enumerate}
		Then the sequence 
		\[f_n:=g_n * 1_{C_{r, \frac{1}{n}}}\qquad \forall n\in\NN\,,\] 
		is an approximate $r$-circle.
	\end{lemma}

	\begin{proof}
	By Theorem~\ref{thm:Gra08}, since $g_n$ is a Schwartz function and $1_{C_{r,\frac{1}{n}}}$ is compactly supported, we have that $f_n = g_n*1_{C_{r,\frac{1}{n}}}\in \cS(\RR^2)$ for each $n\in\NN$.  Next, we need to show that the sequence $\seq{f_n}$ satisfies conditions (C1) and (C2) of Definition~\ref{def:r-circle}.  We do will this in three steps. First, let us define: 
    \[
    u_n \deq (g_n \cdot 1_{B_{\frac{1}{n}}})*1_{C_{r,\frac{1}{n}}}\,, \qquad  v_n \deq f_n - u_n = (g_n\cdot 1_{\RR^2\backslash B_{\frac{1}{n}}})*1_{C_{r,\frac{1}{n}}}\,.
    \] 
    Then $f_n = u_n + v_n$.  \medskip
    
     \noindent \underline{\textit{Step 1.}} We show that there exists a constant $C >0$  and some sequence $a_n \to 0$ such that for all $n\in\NN$ and $x\in\RR^2$ we have 
    \begin{align*}
       0 \leq u_n(x) \leq C\,, \qquad  0 \leq v_n(x) \leq a_n \,.
    \end{align*}    
Observe that:
		\[
		u_n(x) = \int_{\RR^2} g_n(x-y)1_{B_{\frac{1}{n}}}(x-y) 1_{C_{r,\frac{1}{n}}}(y)\,\dd y = \int_{x-C_{r,\frac{1}{n}}} g_n(y) 1_{B_{\frac{1}{n}}}(y)\,\dd y\,,
		\]
		where
		\[
		x-C_{r,\frac{1}{n}}= \lt\{ y \in \RR^2 \,:\, r-\frac{1}{n} < \|x-y \| < r+\frac{1}{n}\rt\}\,.
		\]
		From this and our assumptions we get
		\[
		0 \stackrel{(i)}{\leq} u_n(x) \leq \int_{B_{\frac{1}{n}}} g_n(y)\,\dd y \stackrel{(ii)}{\leq} C < \infty \qquad \forall x\in\RR^2\,,
		\]
		where the last inequality uses that $\int_{B_{\frac{1}{n}}} g_n(x)\,\dd x$ is a convergent sequence, and hence is bounded by some $C>0$.  
	
        Next, by assumptions (i) and (iii) we have:
		\[
		0\stackrel{(i)}{\leq} v_n(x) \stackrel{(iii)}{\leq} \|\varphi\|_\infty \cdot \mbox{\text{Area}}(C_{r,\frac{1}{n}})=: a_n \rightarrow 0\,.
		\]

	\noindent \underline{\textit{Step 2.}} We show that $u_n(x) \rightarrow 1_{C_r}(x)$ for all $x\in\RR^2$.   To prove this, we split the problem into two cases.
	\begin{itemize}
		\item \emph{Case 1:} $\|x\|\neq r$. \smallskip
		
		If $\|x\|\neq r$, then there exists an $N\in\NN$ such that $\overline{B_{\frac{1}{n}}}\cap (x-C_{r,\frac{1}{n}}) = \emptyset$ for all $n\geq N$. 
        From this and the observation that $g_n 1_{B_{\frac{1}{n}}}$ is supported in $\overline{B_{\frac{1}{n}}}$, we get
		\[
		u_n(x) = \int_{x-C_{r,\frac{1}{n}}} g_n(y)1_{B_{\frac{1}{n}}}(y)\,\dd y = \int_{x-C_{r,\frac{1}{n}}} 0 \,\dd y = 0\,,
		\]
		for all $n\geq N$, so $u_n(x)\rightarrow 0$. 
		
		\item \emph{Case 2:} $\|x\|= r$.\smallskip
		
		If $\|x\| = r$, then we have $B_{\frac{1}{n}}\subseteq x - C_{r,\frac{1}{n}}$ for every $n\in\NN$.  From this, we see that
		\[
		u_n(x) = \int_{x-C_{r,\frac{1}{n}}} g_n(y)1_{B_{\frac{1}{n}}}(y)\,\dd y = \int_{B_{\frac{1}{n}}} g_n(y)\,\dd y \,,
		\]
		so by (ii), we get $u_n(x)\rightarrow 1$. 
	\end{itemize}
	Moreover, by Step 1, we have that $\lim_{n\rightarrow\infty} v_n(x) = 0$ for all $x$, so it follows that condition (C1) is satisfied. \medskip

	\noindent \underline{\textit{Step 3.}} We prove (C2).  Observe that $u_n$ is compactly supported with
	\[
	\supp(u_n)\subseteq \supp(g_n \cdot 1_{B_{\frac{1}{n}}}) + \supp(1_{C_{r,\frac{1}{n}}}) \subseteq \overline{B_{\frac{1}{n}}} + \overline{B_{r + \frac{1}{n}}} \subseteq \overline{B_{r + 2}} \qquad \forall n\in\NN\,.
	\]
	Fix some $\psi \in \Cc^\infty(\RR^2)$ such that $\psi \geq 1_{\overline{B_{r+2}}}$. Then, using positivity, we have 
	\begin{equation*}
		0 \stackrel{(i)}{\leq} u_n \leq C 1_{\overline{B_{r+2}}} \leq C \psi \in\cS(\RR^2) \qquad \forall n\in\NN\,,
	\end{equation*}
	and 
	\begin{align*}
		0 \stackrel{(i)}{\leq} v_n  \leq \varphi*\psi \in\cS(\RR^2)\,.
	\end{align*}
	It follows that 
	\[
	0 \leq f_n = u_n + v_n  \leq C \psi+ \varphi*\psi \in \cS(\RR^2) \qquad \forall n\in\NN\,,
	\]
 Therefore, (C2) is satisfied.
\end{proof}

\begin{remark}
    The same proof shows that Lemma~\ref{lem:alt-rcircle} also holds in $\RR^d$.
\end{remark}

Next, we would like to use this to derive a result similar to Theorem~\ref{thm:diffraction-formula} that gives a formula for $\hatmu(C_r)$ for all $\mu\in\STM(\RR^2)$.  We will need the Fourier transform of the characteristic function $1_{B_R}$.  It turns out that this can be expressed in terms of the following Bessel function $J_1$:

\begin{definition}
    For each $t\geq 0$, define
    \[
    J_1(t) \deq \begin{cases}
        \frac{1}{t}\int_0^t s J_0(s) \,\dd s & t > 0 \\
        0 & t = 0\,.
    \end{cases}
    \]
    This function is the \emph{Bessel function of the first kind and order 1}.  
\end{definition}

The above definition of $J_1$ is equivalent to the integral formula
\[
\int_{0}^{x} u J_0(u) \dd u = x J_1(x) \qquad \forall x \in \RR \,, 
\]
which we will use below.

\begin{remark}
$J_1$ is usually defined as
    \[
    J_1(t) = \frac{1}{2\pi}\int_{-\pi}^\pi e^{i(\theta - t\sin\theta)}\,\dd \theta\,.
    \]
The equivalence of these definitions is most easily seen from the well-known series representations for $J_0$ and $J_1$ \cite[p.\,45]{Wat1944}.
\end{remark}

Note that, since $1_{B_R}$ has compact support, its Fourier transform $\reallywidehat{1_{B_R}}$ must be a continuous function.  We get the following formula:

\begin{lemma}\label{lem:disk}
For any $R>0$, $\reallywidehat{1_{B_R}}$ is the continuous function given by
\begin{equation}
\reallywidehat{1_{B_R}}(y) = 
\begin{cases}
    \frac{R}{\|y\|}J_1(2\pi R\|y\|) & y \neq 0 \\
    \pi R^2 & y = 0
\end{cases} \qquad \forall y\in\RR^2\,.
\end{equation} 
\end{lemma}

\begin{proof}
Let $R>0$ and $y\in\RR^2$.  When $y=0$, we have
$$
\reallywidehat{1_{B_R}}(0) = \int_{\RR^2}1_{B_R}(x)\,\dd x = \pi R^2\,.
$$
%Next, when $y\neq 0$, converting to polar coordinates, we get: 
%\begin{align*}
%    \reallywidehat{1_{B_R}}(y) &= \int_{\RR^2} 1_{B_R}(x) e^{-2\pi i x\cdot y}\,\dd x = \int_0^R t  \int_{-\pi}^\pi e^{-2\pi it\|y\|\cos(\theta)}\,\dd\theta \dd t \\
%    &= \int_0^R 2\pi t J_0(2\pi t\|y\|) \dd t \stackrel{u = 2\pi t \|y\|}{=\joinrel=\joinrel=\joinrel=\joinrel=\joinrel=} \frac{1}{2\pi\|y\|^2} \int_0^{2\pi R\|y\|} u J_0(u)\,\dd u\,.
%    \end{align*}
Next consider the case when $y \neq 0$.  
%The derivative in \eqref{eq:J_1-derivative} gives us the following:
%\[
%\int_{0}^{2\pi R \|y\|} u J_0(u)\,\dd u = u J_1(u) \big|_0^{2\pi R \|y\|} = 2\pi R \|y\|J_1(2\pi R \|y\|)\,.
%\]
Since $1_{B_R}$ is rotational invariant we get
\begin{align*}
 \reallywidehat{1_{B_R}}(y) &=  \reallywidehat{1_{B_R}}(\| y \|,0 ) = \int_{B_R} e^{-2 \pi i x \cdot (\| y \|,0 )} \dd x  = \int_{0}^R \int_{-\pi}^\pi e^{-2 \pi i r \|y \| \cos(\theta)} r \dd \theta \dd r \\
 &= 2 \pi  \int_0^R r J_0(2 \pi r \|y \|) \dd r \stackrel{u=2 \pi r \|y \|}{=\joinrel=\joinrel=\joinrel=\joinrel=\joinrel=\joinrel=\joinrel=}\frac{1}{2 \pi \| y \|^2} \int_0^{2 \pi R \|y \|} u J_0(u) \dd u =\frac{R}{\| y \|} J_1(2 \pi R \|y \|) \,,
\end{align*}
as claimed.  

\end{proof}

Next, observe that whenever $n>\frac{1}{r}$, we have
$$
\reallywidehat{1_{C_{r,\frac{1}{n}}}} = \reallywidehat{1_{B_{r+\frac{1}{n}}}} - \reallywidehat{1_{B_{r-\frac{1}{n}}}}\,.
$$
This lets us evaluate $\reallywidehat{1_{C_{r,\frac{1}{n}}}}$ using Lemma~\ref{lem:disk}:

\begin{corollary}\label{cor:annulus-FT}
For any $r>0$ and $n > \frac{1}{r}$, we have 
\begin{equation}\label{eq:annulus-FT}
\reallywidehat{1_{C_{r,\frac{1}{n}}}}(y)  = 
\begin{cases}
    \frac{nr+1}{n\|y\|}J_1\lt(\frac{2\pi \|y\|(nr+1)}{n}\rt) - \frac{nr-1}{n\|y\|}J_1\lt(\frac{2\pi \|y\|(nr-1)}{n}\rt) & y\neq 0 \\
    \pi \lt(r+\frac{1}{n}\rt)^2 - \pi\lt(r-\frac{1}{n}\rt)^2 & y = 0 
\end{cases} \qquad \forall y\in\RR^2\,.
\end{equation}\qed
\end{corollary}

As shown above, we have an explicit formula for the Fourier transform $\reallywidehat{1_{C_{r,\frac{1}{n}}}}$.  Thus, given a Schwartz function $g_n$ as in Lemma~\ref{lem:alt-rcircle}, we can use the convolution theorem to get
$$\reallywidehat{g_n*1_{C_{r,\frac{1}{n}}}} = \reallywidehat{g_n}\cdot \reallywidehat{1_{C_{r,\frac{1}{n}}}}\,,$$ 
which would allow us to compute explicitly the Fourier transform $\reallywidehat{g_n*1_{C_{r,\frac{1}{n}}}}$ provided that $\reallywidehat{g_n}$ can also be computed explicitly. In particular, combining Lemma~\ref{lem:DCT}, Lemma~\ref{lem:alt-rcircle}, and Corollary~\ref{cor:annulus-FT} gives us the following general theorem:

\begin{theorem}\label{thm:alt-formula-general}
    Let $r>0$. If $\seq{g_n}$ is a sequence in $\cS(\RR^2)$ satisfying the conditions of Lemma~\ref{lem:alt-rcircle} and $\mu\in\STM(\RR^2)$, then we have
    \begin{equation}\label{eq:alt-formula-general}
    \hatmu(C_r) = \lim_{n\rightarrow\infty} \int_{\RR^2}\reallywidehat{g_n}(y) \reallywidehat{1_{C_{r,\frac{1}{n}}}}(y) \,\dd \mu(y)\,,
    \end{equation}
    where $\reallywidehat{1_{C_{r,\frac{1}{n}}}}(y)$ is given by \eqref{eq:annulus-FT}. \qed
\end{theorem}

\begin{remark}
    The above formula \eqref{eq:alt-formula-general} also holds when $\seq{g_n}$ is a sequence of compactly supported functions in $L^1(\RR^2)\cap L^2(\RR^2)$ satisfying
    $$
    \supp(g_n)\subseteq \overline{B_{\frac{1}{n}}} \qquad \forall n\in\NN\,,
    $$
    and conditions (i) and (ii) of Lemma~\ref{lem:alt-rcircle}, provided that $\hatmu$ is translation bounded.  This can be seen from \cite[Theorem~3.12]{RS2017a} and the Fourier transformability of $\mu$ as a measure (see Remark~\ref{r1}).
\end{remark}

We conclude with an example using a sequence of Gaussian functions and Theorem~\ref{thm:alt-formula-general} to produce an alternate formula \eqref{eq:alt-diffraction-formula} for the intensity of a circle.

\begin{example}
    Consider the sequence of Gaussian functions 
    $$
    g_n(x) \deq \frac{n^4}{\pi}e^{-n^4 \|x\|^2} \qquad \forall x\in\RR^2, \ \forall n\in\NN\,.
    $$       
    Then $g_n\geq 0$ and 
    \begin{align*}
        \int_{B_{\frac{1}{n}}} g_n(x)\,\dd x &=   \frac{n^4}{\pi}\int_{\| x\| \leq \frac{1}{n}} e^{-n^4 \|x\|^2}\,\dd x \\
        &\stackrel{y=n^2x}{=\joinrel=\joinrel=\joinrel=\joinrel=} \frac{1}{\pi} \int_{\| y\| \leq n}  e^{-\|y\|^2}\,\dd y  \rightarrow \frac{1}{\pi} \int_{\RR^2}e^{-\|y\|^2}\,\dd y  = 1\,,
    \end{align*}
    so $\{g_n\}_{n\in\NN}$ satisfies conditions (i) and (ii) of Lemma~\ref{lem:alt-rcircle}.

    Next, to verify condition (iii), we demonstrate that, for the Schwartz function 
    \[
    \varphi(x) \deq \frac{4}{\pi} e^{1-\|x\|^2} \qquad \forall x \in \RR^2\,,
    \]
    we have $g_n(x) \leq \varphi(x)$ whenever $x\in \RR^2\backslash B_{\frac{1}{n}}$. Equivalently, we want to verify that
    \begin{equation}\label{eq:Gaussian-bound}
    \frac{n^4}{4} \leq e^{1+(n^4 - 1)t^2} \qquad \forall t \geq \frac{1}{n}\,.
    \end{equation}
    Indeed, for $t \geq \frac{1}{n}$, we have
    \[
    1+(n^4 - 1)t^2 \geq 1 + n^2 - \frac{1}{n^2} \geq n^2\,,
    \]
    so
    \[
    \ln\lt(\frac{n^4}{4}\rt) = 2\ln\lt(\frac{n^2}{2}\rt) \leq 2 \cdot \frac{n^2}{2} = n^2 \leq 1+(n^4 - 1)t^2 = \ln \lt( e^{1+(n^4 - 1)t^2} \rt)\,.
    \]
    This verifies \eqref{eq:Gaussian-bound}.  

    Thus, we have shown that $\seq{g_n}$ satisfies the conditions of Lemma~\ref{lem:alt-rcircle}, so given a measure $\mu\in\STM(\RR^2)$ and some $r>0$, we can apply Theorem~\ref{thm:alt-formula-general} to obtain the following formula for $\hatmu(C_r)$:
    \begin{multline}\label{eq:alt-diffraction-formula}
    \hatgam(C_r)  
    =  \lim_{n\rightarrow\infty} \int_{\RR^2} e^{-\pi^2 \|y\|^2/n^4} \\
    \times \lt[\frac{nr+1}{n\|y\|}J_1\lt(\frac{2\pi \|y\|(nr+1)}{n}\rt) - \frac{nr-1}{n\|y\|}J_1\lt(\frac{2\pi \|y\|(nr-1)}{n}\rt)\rt]
    \,\dd \gamma(y)\,.
    \end{multline}
    for any $\mu\in\STM(\RR^2)$ and any $r>0$.
\end{example}

\section*{Acknowledgements}
E.R.K. was supported by NSERC via the NSERC CGS D scholarship. N.S. was supported by the NSERC Discovery grants 2020-00038 and 2024-04853. The authors are grateful to Michael Baake and Noel Murasko
%and Shigeki Akiyama 
for some discussions and suggestions which improved the quality of the manuscript.

\printbibliography

@article{AG1974,
  author = {Argabright, L. N. and Gil de Lamadrid, J.},
  title = {{F}ourier analysis of unbounded measures on locally compact abelian groups},
  year = {1974},
  journal = {{Me}moirs of the {A}merican {M}athematical {S}ociety},
  number = {145},
  address = {AMS, Providence, RI}
}

@article{BFG2007,
  title = {A Radial Analogue of {{Poisson}}'s Summation Formula with Applications to Powder Diffraction and Pinwheel Patterns},
  author = {Baake, Michael and Frettl{\"o}h, Dirk and Grimm, Uwe},
  year = {2007},
  journal = {Journal of Geometry and Physics},
  volume = {57},
  number = {5},
  pages = {1331--1343}
}

@article{BFG2007a,
  title = {Pinwheel Patterns and Powder Diffraction},
  author = {Baake, Michael and Frettl{\"o}h, Dirk and Grimm, Uwe},
  year = {2007},
  journal = {Philosophical Magazine},
  volume = {87},
  number = {18-21},
  eprint = {https://doi.org/10.1080/14786430601057953},
  pages = {2831--2838},
  publisher = {{Taylor \& Francis}},
  doi = {10.1080/14786430601057953}
}

@article{BG2003,
  title = {A note on shelling},
  author = {Baake, Michael and Grimm, Uwe},
  year = {2003},
  journal = {Discrete \& Computational Geometry},
  volume = {30},
  pages = {573-589}
}

@book{BG2013,
  title = {{Aperiodic Order.  Vol. 1:  A Mathematical Invitation}},
  author = {Baake, Michael and Grimm, Uwe},
  year = {2013},
  publisher = {{Cambridge University Press}},
  address = {{Cambridge}}
}

@book{BG2017,
  title = {{Aperiodic Order.  Vol. 2:  Crystallography and Almost Periodicity}},
  author = {Baake, Michael and Grimm, Uwe},
  year = {2017},
  publisher = {{Cambridge University Press}},
  address = {{Cambridge}}
}

@article{BHS2017,
title = {On weak model sets of extremal density},
journal = {Indagationes Mathematicae},
volume = {28},
number = {1},
pages = {3-31},
year = {2017},
note = {Special Issue on Automatic Sequences, Number Theory, and Aperiodic Order},
author = {Michael Baake and Christian Huck and Nicolae Strungaru}
}

@book{BM2000,
	editor = {Michael Baake and Robert V. Moody},
	title = {Directions in Mathematical Quasicrystals},
	year = {2000},
	series = {CRM Monograph Series},
	volume = {13},
	publisher = {AMS},
	address = {Providence, RI}
}

@article{BM2004,
year = {2004},
author = {Michael Baake and Robert V. Moody},
title={Weighted {D}irac combs with pure point diffraction},
journal = {Journal für die reine und angewandte Mathematik (Crelle's Journal)},
volume = {573},
pages = {61-94}
}

@book{BO1999,
Author = {Bender, Carl M. and Orszag, Steven A.},
Publisher = {Springer New York},
Title = {Advanced Mathematical Methods for Scientists and Engineers I : Asymptotic Methods and Perturbation Theory},
address = {New York},
Year = {1999}
}

@article{BS2022,
      title={A note on tempered measures}, 
      author={Michael Baake and Nicolae Strungaru},
      year={2023},
      pages={15-30},
	journal={Colloquium Mathematicum},
	number={172}
}

@inbook{MS2017,
  title = {Almost Periodic Measures and Their {F}ourier Transforms},
  author = {Moody, Robert V. and Strungaru, Nicolae},
  chapter = 4,
  crossref="BG2017"
}

@book{Ped1995,
  title = {Analysis Now},
  author = {Pedersen, Gert Kjaerg{\aa}rd},
  year = {1995},
  series = {Graduate Texts in Mathematics},
  edition = {Rev. print., corr. 2. print},
  number = {118},
  publisher = {{Springer}},
  address = {{New York Berlin Heidelberg}},
  isbn = {978-1-4612-1007-8 978-0-387-96788-2 978-3-540-96788-0},
  langid = {english}
}

@article{RS2017,
  title = {A Short Guide to Pure Point Diffraction in Cut-and-Project Sets},
  author = {Richard, Christoph and Strungaru, Nicolae},
  year = {2017},
  journal = {Journal of Physics A: Mathematical and Theoretical},
  volume = {50},
  number = {15},
  pages = {154003}
}

@article{RS2017a,
  title = {Pure Point Diffraction and {{Poisson}} Summation},
  author = {Richard, Christoph and Strungaru, Nicolae},
  year = {2017},
  journal = {Annales Henri Poincar\'e},
  volume = {18},
  number = {12},
  primaryclass = {math-ph},
  pages = {3903--3931}
}

@inbook{Str2017,
	author = {Nicolae Strungaru},
	title = {Almost periodic pure point measures},
	chapter = 5,
	crossref = "BG2017"
}

@article{SS2021,
  title = {On the (Dis)Continuity of the {{Fourier}} Transform of Measures},
  author = {Spindeler, Timo and Strungaru, Nicolae},
  year = {2021},
  journal = {Journal of Mathematical Analysis and Applications},
  volume = {499},
  number = {2},
  pages = {125062}
}

@article{BGJR2000,
title = {Averaged shelling for quasicrystals},
journal = {Materials Science and Engineering: A},
volume = {294-296},
pages = {441-445},
year = {2000},
author = {Michael Baake and Uwe Grimm and Dieter Joseph and Przemysław Repetowicz}
}

@online{FGH,
author = {Dirk Frettlöh and Edmund Harriss and Franz Gähler},
title = {Tilings Encyclopedia},
url= {https://tilings.math.uni-bielefeld.de/},
date=2024
}

@article{Fre2008,
  title = {Substitution Tilings with Statistical Circular Symmetry},
  author = {Frettl{\"o}h, Dirk},
  year = {2008},
  journal = {European Journal of Combinatorics},
  volume = {29},
  number = {8},
  pages = {1881--1893}
}

@article{GD2011,
  title = {Some Comments on Pinwheel Tilings and Their Diffraction},
  author = {Grimm, Uwe and Deng, Xinghua},
  year = {2011},
  journal = {Journal of Physics: Conference Series},
  volume = {284},
  number = {1},
  pages = {012032}
}

@book{Gra2008,
  title = {Classical {{Fourier}} Analysis},
  author = {Grafakos, Loukas},
  year = {2008},
  series = {Graduate Texts in Mathematics},
  edition = {2nd ed},
  number = {249},
  publisher = {{Springer}},
  address = {{New York}}
}

@book{GR2000,
	author = {Gradshteyn, Izrail Solomonovich and Ryzhik, Iosif Moiseevich},
	title = {Table of Integrals, Series, and Products},
	editor = {Jeffrey, Alan and Zwillinger, Daniel},
	year = {2000},
	publisher = {Academic Press},
	address = {San Diego},
	edition = 6
}

@article{Hof1995,
  title = {On Diffraction by Aperiodic Structures},
  author = {Hof, Albertus},
  year = {1995},
  journal = {Communications in Mathematical Physics},
  volume = {169},
  number = {1},
  pages = {25--43},
  publisher = {{Springer-Verlag}},
  address = {{Berlin}}
}

@article{KR2016, 
title={Dynamics on the graph of the torus parametrization}, volume={38}, 
number={3}, 
journal={Ergodic Theory and Dynamical Systems}, 
author={Keller, Gerhard and Richard, Christoph}, 
year={2016}, 
pages={1048–1085}}

@misc{LSS2020,
      title={Pure point diffraction and mean, {B}esicovitch and {W}eyl almost periodicity}, 
      author={Daniel Lenz and Timo Spindeler and Nicolae Strungaru},
      year={2020},
      howpublished={arXiv:2006.10821}
}

@article{LSS2024, 
title={Pure point spectrum for dynamical systems and mean, {B}esicovitch and {W}eyl almost periodicity},
volume={44}, 
number={2}, 
journal={Ergodic Theory and Dynamical Systems}, 
author={Lenz, Daniel and Spindeler, Timo and Strungaru, Nicolae}, 
year={2024}, 
pages={524–568}}

@article{Rad1994,
  title = {The Pinwheel Tilings of the Plane},
  author = {Radin, Charles},
  year = {1994},
  journal = {Annals of Mathematics},
  volume = {139},
  number = {3},
  eprint = {2118575},
  eprinttype = {jstor},
  pages = {661--702},
  publisher = {{Annals of Mathematics}}
}

@article{Sad1998,
  title = {Some generalizations of the pinwheel tiling},
  author = {Sadun, Lorenzo},
  year = {1998},
  journal = {Discrete \& Computational Geometry},
  volume = {20},
  pages = {79--110}
}

@inbook{Sch2000,
	author = {Martin Schlottmann},
	title = {Generalized Model Sets and Dynamical Systems},
	pages = {144-147},
	crossref = "BM2000"
}

@article{SBGC1984,
  title = {Metallic Phase with Long-Range Orientational Order and No Translational Symmetry},
  author = {Shechtman, D. and Blech, I. and Gratias, D. and Cahn, J. W.},
  journal = {Physical Review Letters},
  volume = {53},
  number= {20},
  pages = {1951--1953},
  numpages = {0},
  year = {1984},
  publisher = {American Physical Society}
}

@article{Str2005,
	author = {Nicolae Strungaru},
	title = {Almost periodic measures and long-range order in {M}eyer sets},
	journal = {Discrete {\&} Computational Geometry},
	volume = {33},
	pages = {483-505},
	year = {2005}
}

@article{Str2020,
	author = {Nicolae Strungaru},
	title = {On the {F}ourier analysis of measures with {M}eyer set support},
	journal = {Journal of Functional Analysis},
	volume = {278},
	number = {6},
	pages = {108404},
	year = {2020}
}

@book{Wat1944,
Author = {Watson, George Neville},
Publisher = {Cambridge University Press},
Title = {A treatise on the theory of Bessel functions},
Year = {1944},
Address = {Cambridge}
}

@article{MPS2006,
  title = {Circular symmetry of pinwheel diffraction},
  author = {Moody, Robert V. and Postnikoff, Derek and Strungaru, Nicolae},
  year = {2006},
  journal = {Annales Henri Poincar\'e},
  volume = {7},
  pages = {711--730}
}

\end{document}